\documentclass[10pt]{article}
\usepackage{fullpage} %
\usepackage{url} %
\usepackage{amsmath} %
\usepackage[backrefs,lite,abbrev]{amsrefs} %
\usepackage{color} %
\usepackage{tikz} %
\usepackage{amssymb,latexsym} %
\usepackage{framed}
\definecolor{shadecolor}{rgb}{.9,.9,.9}%{\color{yellow}}
\newcommand{\betabox}{{\gamma}}
\newcommand{\junk}[1]%{}
{\begin{shaded}\slshape  #1 \end{shaded}}
\newcommand{\lsty}[1]{{\large\bf #1}}

\newtheorem{theorem}{Theorem}[section]
\newtheorem{lemma}[theorem]{Lemma}
\newtheorem{corollary}[theorem]{Corollary}
\newtheorem{definition}[theorem]{Definition}
\newtheorem{proposition}[theorem]{Proposition}
\newtheorem{example}[theorem]{Example}
\newtheorem{rem}[theorem]{Remark}
\newtheorem{notation}[theorem]{Notation}
\newcommand {\eop}        {\hfill $\Box$}
\newenvironment{proof}%
{ \medskip
  \rmfamily \noindent
  {\bf Proof.\/} } {\eop \par \medskip}

\newcommand{\B} {\ensuremath{\mathcal{B}}}
\newcommand{\R}{\ensuremath{\mathbb{R}}}

\newcommand{\bx}{\mathrm{box}}
\newcommand{\G}{\mathcal{G}}
\newcommand{\GG}{\widehat{G}}

\let\sub=\subseteq

\let\bs=\backslash

\begin{document}
\author{Michael Abrahams
\thanks{Vassar '10}, 
Meg Lippincott
\thanks{Vassar '09},
and Thierry Zell
\thanks{The Donald and Helen Schort school of Mathematics and 
Computing Sciences,  Lenoir-Rhyne University, Hickory NC 28603}}

\date{August 23, 2009}
\title{On $(2,3)$-agreeable Box Societies}

\maketitle

\begin{abstract}
The notion of $(k,m)$-agreeable society was introduced by Deborah 
Berg et al.:% in~\cite{agree}: 
a family of convex subsets of $\R^d$ is 
called $(k,m)$-agreeable if any subfamily of size $m$ contains at least one non-empty $k$-fold intersection. 
%In~\cite{agree}, 
In that paper,
the $(k,m)$-agreeability of a convex family was shown to imply the existence of a subfamily of size 
$\beta n$ with non-empty intersection, where $n$ is the size of the original family and $\beta\in[0,1]$ is an explicit constant depending only on $k,m$ and $d$. The quantity $\beta(k,m,d)$ is called the minimal \emph{agreement proportion} for a $(k,m)$-agreeable family in $\R^d$.

If we only assume that the sets are convex,
simple examples show that 
$\beta=0$ for $(k,m)$-agreeable families in $\R^d$ where $k<d$. 
In this paper, we introduce new techniques to find positive lower bounds  when restricting our attention to  families 
of $d$-boxes, i.e. cuboids with sides parallel 
to the coordinates hyperplanes. We  derive explicit formulas
for the first non-trivial case: the case of $(2,3)$-agreeable families of $d$-boxes with $d\geq 2$.
\end{abstract}

%\tableofcontents

\section{Introduction}

The article~\cite{agree} introduced the concept of geometric approval
voting, where a \emph{platform} is a point in $\R^d$ and a 
\emph{vote} can be any
convex subset, representing all the platforms deemed acceptable by
that particular voter. 
(The convexity assumption is a way to require our voters to be reasonable: the fact that all votes  
contain every point on a segment with both endpoints in the vote
means that any platform obtained  as a compromise between two acceptable positions is again deemed acceptable.)
The main question addressed in~\cite{agree} was,
given a collection of votes, to
find the largest number of overlapping votes, and thus the largest
number of voters that could be satisfied by the adoption of any single 
platform.

More specifically, the authors concentrated on what they termed
\emph{$(k,m)$-agreeable societies}, where any group of $m$ voters
contains $k$ or more  who can agree on a common platform. 
Their main goal was to obtain lower bounds on the \emph{agreement
proportion} (the ratio number of satisfied voters over total number of voters)
in terms of $k,m$ and $d$ only.
Using the version of the fractional Helly theorem due to Kalai~\cite{kalai},
they showed that if the society contains $n \geq m$ votes, all of which are
convex subsets of $\R^d$, then there exsits a platform contained in at
least $\beta(k,m,d)\ n$ votes, where the \emph{proportion} 
$\beta(k,m,d)$ verifies:
\begin{equation}\label{eq:cvxbeta}
\beta(k,m,d) \geq 1-
\left[1-\frac{\binom{k}{d+1}}{\binom{m}{d+1}} \right]^{\frac{1}{d+1}}.
\end{equation}
Given that the fractional Helly theorem cannot use information on the
number of $k$-fold intersection when $k\leq d$, it is no surprise that
this lower bound is positive only when $k \geq d+1$.

If the general convex case requires detailed information about the
whole \emph{nerve complex} of the arrangement of votes, the \emph{intersection graph} does capture the complexity of the
whole arrangement in the special case when the votes are \emph{boxes},
i.e. parallelotopes whose sides are parallel to the coordinate
axes. This case was also addressed in~\cite{agree}, and purely
graph-theoretic considerations yielded a sharp bound
of $k/m$ for the agreement proportion in the 
\emph{strong agreement case:} the situation of $(k,m)$-agreeability where $m \leq 2k-2$.
(The result proved in~\cite{agree} for this case $m \leq 2k-2$ 
is in fact substantially stronger: if the number of boxes is $n$, there is an overlap of at least
$n-m+k$ boxes, so the actual agreement proportion \emph{starts} at $k/m$ and \emph{increases} to 1 with the number $n$ of boxes.)

The case of societies of $(2,m)$-agreeable $d$-boxes does not
fall in the strong agreement category, and it
is left essentially open in~\cite{agree}.
In fact, it is not
even clear at the outset 
that there is a positive agreement proportion for 
$(2,m)$-agreeable $d$-box arrangements when
$m \geq 3$ and $d\geq2$, since the lower bound given 
by~\eqref{eq:cvxbeta} is zero in that case.
In this paper, we tackle the $(2,3)$-agreeable case and we prove the
following result.

\begin{theorem}\label{thm:main}
For any $d \geq 1$, 
any $(2,3)$-agreeable $d$-box society  has an 
agreement proportion of at least $(2d)^{-1}$.
\end{theorem}

The remainder of the paper is organized as follows.
\begin{description}
\item[Linear Case.] The material in Section~2 is independent from the rest of the paper: it presents an elementary proof of the fact that $(2,3)$-agreeable arrangements of interval have agreement proportion $1/2$. 
\item[Preliminaries.] Section~3 introduces basic notations and definitions regarding arrangements of boxes and their intersection graphs.
\item[Degree Bounds.] Section~4 establishes lower- and upper-bounds 
on the degrees of vertices of $(2,3)$-agreeable graphs with bounded clique number. A classification of the small cases is given, and 
we prove that positive lower bounds do exist for all $d$.
\item[Main Result.] In Section~5, we establish the specific values of the lower bound stated in Theorem~\ref{thm:main}. The proof uses a lower bound on boxicity taken from Adiga et al.~\cite{lower}.
Section~6 presents a few questions left open by our work.
\item[Appendix.] We finish the paper with an entirely different lower bound proof. The bounds obtained are somewhat weaker, but we believe that the technique, borrowing important ideas about arrangement of boxes from Eckhoff's work~\cite{boxes1}, is interesting in its own right in view of its applicability in other settings.
\end{description}
Throughout the paper, all arrangements of boxes are assumed to be $(2,3)$-agreeable. Many of the definitions and results could easily be 
extended to the $(k,m)$-agreeable case; this level of generality was eschewed in order to keep notations simple and legible. The only step for which $(2,3)$-agreeability is crucial is in establishing the lower bound of Section~4.

\bigskip
\noindent 
{\bf Acknowledgments.}  This paper originated with a research
project that took place during the  2008 
Undergraduate Research Summer Institute at Vassar College, 
where the first two authors were students and the last author was a visiting professor.  
The authors are grateful to the institute for its support, and extend
special thanks to Professor Frank and her own URSI group for helping to foster
a stimulating mathematical environment.

The software \emph{Mathematica}, and especially the \emph{Combinatorica} package, 
proved invaluable in the study of examples for this paper.

\begin{notation}
Throughout this paper, $G$ denotes a simple, undirected graph.
The sets $V(G)$ and $E(G)$ are respectively the sets of vertices and edges of $G$, and we let $n=\#V(G)$. Recall that any subset $W$ of $V(G)$ gives rise to the subgraph $G[W]$ \emph{induced} by $W$, which is the graph which has $W$ as its set of vertices, and has for edges all the edges of $E(G)$ with both endpoints in $W$.

A \emph{clique} in $G$ is any subset of $V(G)$ that induces a complete subgraph, and the size of the largest clique is called the \emph{clique number} of $G$ and denoted by $\omega(G)$.
\end{notation}

\section{The Linear Case}
The intersection graphs associated to arrangements of intervals in the
line are \emph{perfect graphs}. 
This allowed the authors of~\cite{agree} to prove the non-trivial fact: for any $(k,m)$-agreeable arrangement of intervals, the agreement
number is at least $(n-R)/Q$, where $Q$ and $R$ denote respectively
the quotient and the remainder of the euclidean division of $m-1$ by $k-1$. This
lower bound is sharp, and it implies that any $(k,m)$-agreeable collection of intervals must have an agreement proportion
$$\beta(k,m,1) \geq\frac{k-1}{m-1}.$$  
In particular, the above implies that 
any $(2,3)$-agreeable collection of intervals
has agreement proportion at least $1/2$. 
This substantially improves the general case bound given 
in the formula~\eqref{eq:cvxbeta}, 
which for $d=1$ in the $(2,3)$-agreeable 
setting yields an agreement  proportion of 
$$ 1-\sqrt{\frac{2}{3}} \approx 0.1835. $$
We reprove the bound of $1/2$ using only elementary means.
First, we need to know when the agreement proportion equals~1.

\begin{lemma}
A linear society has agreement proportion~1 if and only if every 
pair of votes intersects.
In the terminology of~\cite{agree}, such an arrangement is called \emph{super-agreeable}.
\end{lemma}

\begin{proof}
This is a special case of Helly's theorem~\cite{helly}, which 
states that for any arrangement of convex sets in $\R^d$, the sets have a non-empty intersection if and only if all $(d+1)$-fold intersections are non-empty.
\end{proof}

\begin{theorem}[\cite{agree}*{Theorem~1}]\label{thm:d1}
The agreement proportion of a linear $(2,3)$-agreeable society is $1/2$. 
\end{theorem}

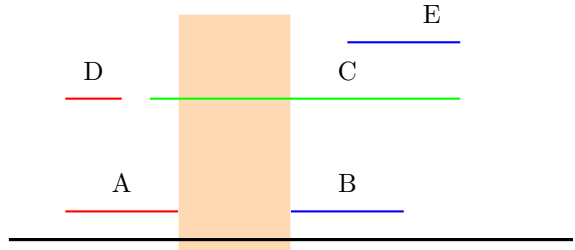
\begin{figure}[ht]
\centering

\begin{tikzpicture}[scale=1.5]
\draw[color=white, fill=orange, fill opacity=0.3] (1.5,-0.1) rectangle (2.5,2); 
\draw[style=very thick] (0,0) -- (5,0);
\draw[color=red, style=thick] (.5,.25) -- (1.5,.25);
\draw[color=blue, style=thick] (2.5, .25) -- (3.5,.25);
\draw[color=green, style=thick] (1.25,1.25) -- (4,1.25);
\draw[color=red, style=thick] (.5,1.25) -- (1,1.25);
\draw[color=blue, style=thick] (3,1.75) -- (4,1.75);
\path (1,.5) node (a) {A};
\path (3,.5) node (b) {B};
\path (3,1.5) node (c) {C};
\path (.75,1.5) node (d) {D};
\path (3.75,2) node (e) {E};
\end{tikzpicture}
\caption{Alice, Bob, and their friends. The shaded area in the middle is shared by all the friends of Bob and Alice, such as $C$.}
\label{fig:d1}
\end{figure}

\begin{proof}
If every pair of votes intersects, Helly's
theorem for intervals implies that the agreement proportion is 1.
So, without loss of generality, we can assume that in our 
one-dimensional
$(2,3)$-agreeable society, there are two non-intersecting intervals
$A$ (Alice's vote) and $B$ (Bob's vote), with $A$ to the left of $B$.

\medskip

The remaining voters can be divided into three categories: those who
only agree with Alice, those who only agree with Bob, and those who
can agree with both Alice and Bob. (There are no voters who agree with neither since that would violate $(2,3)$-agreeability.)
These three categories of voters -- call them friends of Alice, friends of Bob and friends of both --  form super-agreeable groups, where all voters can agree pairwise and thus, by Helly's theorem, all the votes in each group overlap.
Indeed, friends of Alice must
agree with each other, because if two of them did not agree, then
taken together with Bob, we would have three votes containing no
intersecting pair, violating the condition of $(2,3)$-agreeability.
Similarly, voters who only agree with Bob must also agree with each
other. As for votes which overlap with both Alice and Bob's vote, 
they all meet in the interval 
$[\max(A), \min(B)]$ between $A$ and $B$ 
(Figure~\ref{fig:d1}). 
If one of the three categories is empty, we have two super-agreeable groups, one of which must account for at least one half of the voters, and the result holds. 

\medskip

Suppose all three categories are non empty, and let $C$ be a vote containing $[\max(A), \min(B)]$, $D$ be a vote intersecting $A$ only and $E$ be a vote intersecting $B$ only. The three votes must share at least one intersection to respect the $(2,3)$-agreeable condition; and note that if $D\cap E\neq \varnothing$, it implies that the two intersctions with $C$ are also non-empty (all meet in the middle region). If we can find a vote $D$ from a friend of Alice such that $C\cap D=\varnothing$, then we must have $C\cap E\neq \varnothing$, and, replacing $E$ by any other vote $E'$ intersecting $B$, the same reasoning shows 
that $C\cap E'\neq \varnothing$ too. Thus any vote $C$ bridging the gap between Alice and Bob must either meet all the votes that intersect $A$ or all the votes that intersect $B$.
%
%\medskip
Thus, we can assign those bridging votes to Alice or Bob, since they have to overlap with all of the friends of at least one. We can divide the votes into two super-agreeable groups once again. One of those must account for at least half the voters, proving the result.
\end{proof}

\begin{rem}\label{rem:d1}
This theorem is sharp: any society formed by taking $r$
copies of $A$ and $r$ copies of $B$
is $(2,3)$-agreeable with agreement proportion $1/2$.
\end{rem}

\begin{figure}[ht]
\centering

\begin{tikzpicture}[scale=.5]
\draw[fill=orange, fill opacity=0.3] (1,0) rectangle (4,5); 
\draw[fill=blue, fill opacity=0.3] (0,3) rectangle (2,12); 
\draw[fill=purple, fill opacity=0.3] (3,4) rectangle (6,6); 
\draw[fill=green, fill opacity=0.3] (5,5) rectangle (7,10); 
\draw[fill=red, fill opacity=0.3] (1,7) rectangle (6,9); 
\end{tikzpicture}
\caption{A $(2,3)$-agreeable society of $2$-boxes with agreement
proportion $2/5$.}
\label{fig:Z5}
\end{figure}
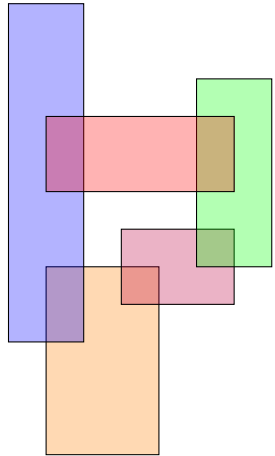

\begin{rem}\label{rem:Z5}
The result of the previous theorem only holds in dimension~1. 
For instance, Figure~\ref{fig:Z5} shows five votes in dimension~2 
arising from $(2,3)$-agreeable voters, yet the agreement proportion is only $2/5$.
\end{rem}

\section{Boxes and Agreeable Graphs}

We introduce some definitions and notations for the two main objects of study: arrangements of boxes and their associated intersection graphs.

\subsection{Arrangements of Boxes and Intersection Graphs}

A \emph{$d$-box} is a subset of $\R^d$ given by the cartesian product
of $d$ closed intervals.
A collection $\B$ of boxes gives rise to a graph in the following fashion.

\begin{definition}
The \emph{intersection graph} $G_\B$ associated to an arrangement
$\B=\{B_1, \dots, B_n\}$ of $d$-boxes is the graph with vertices
$V=\{1, \dots, n\}$ and such that $\{i,j\}$ is an edge if and only if $B_i
\cap B_j \neq \varnothing$.

Conversely, given a simple undirected graph $G$, we can define its
\emph{boxicity} $\bx(G)$: it is the smallest integer $d$ such that
there exists an arrangement of $d$-boxes $\B$ whose intersection graph
is $G$. 
\end{definition}
 Roberts~\cite{roberts} showed that this number is always
finite, and that $\bx(G) \leq \lfloor \#V/2 \rfloor$. 
(Graphs for which this bound is tight are classified in~\cite{trotter}.) 

\begin{rem}\label{rem:subgraph}
By convention, we let $\bx(K_n)=0$ for all $n$ (a $0$-box would be a point). This shows that boxicity does not behave nicely with respect to taking subgraphs. On the other hand, it is clear that boxicity can only decrease when taking \emph{induced subgraphs}, since for any arrangement of $d$-boxes $\B=\{B_1, \dots, B_n\}$ and any 
subset $I \sub \{1, \dots, n\}$, the intersection graph of the sub-arrangement 
$\{B_i\mid i\in I\}$ is simply the graph $G_\B[I]$ induced by 
the vertices $I$ in $G_\B$.
\end{rem}

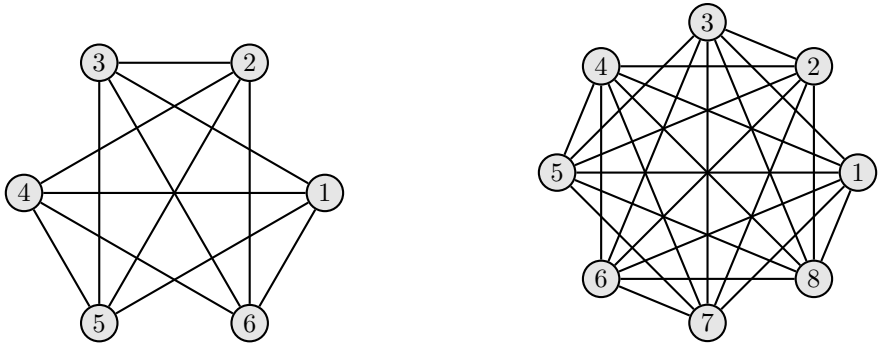
\begin{figure}[ht]
\centering
\begin{tikzpicture}[scale=2, style=thick]
\tikzstyle{every node}=[circle, draw, fill=black!10,
                        inner sep=2pt, minimum width=4pt]
\path (0:1) node (1) {1};
\path (60:1) node (2) {2};
\path (120:1) node (3) {3};
\path (180:1) node (4) {4};
\path (240:1) node (5) {5};
\path (300:1) node (6) {6};
\draw \foreach \x in {3,4,5,6} 
{ (\x) -- (1) };
\draw \foreach \x in {3,4,5,6} 
{ (\x) -- (2) };
\draw \foreach \x in {5,6} 
{ (\x) -- (3) };
\draw \foreach \x in {5,6} 
{ (\x) -- (4) };
\end{tikzpicture}
\hfil
\begin{tikzpicture}[scale=2, style=thick]
\tikzstyle{every node}=[circle, draw, fill=black!10,
                        inner sep=2pt, minimum width=4pt]
\path (0:1) node (1) {1};
\path (45:1) node (2) {2};
\path (90:1) node (3) {3};
\path (135:1) node (4) {4};
\path (180:1) node (5) {5};
\path (225:1) node (6) {6};
\path (270:1) node (7) {7};
\path (315:1) node (8) {8};
\draw \foreach \x in {3,4,5,6,7,8} 
{ (\x) -- (1) };
\draw \foreach \x in {3,4,5,6,7,8} 
{ (\x) -- (2) };
\draw \foreach \x in {5,6,7,8} 
{ (\x) -- (3) };
\draw \foreach \x in {5,6,7,8} 
{ (\x) -- (4) };
\draw \foreach \x in {7,8} 
{ (\x) -- (5) };
\draw \foreach \x in {7,8} 
{ (\x) -- (6) };
\end{tikzpicture}

\caption{The complete partite graphs $K_3(2)$ and $K_4(2)$.}\label{fig:K}
\end{figure}

\begin{example}\label{ex:K}
Note that the
bound $\bx(G) \leq \lfloor \#V/2 \rfloor$ remains sharp,
\emph{even if we restrict our attention to $(2,3)$-agreeable graphs.}
Indeed, for any $d\geq 1$, let $K_d(2)$ be the complete $d$-partite graph on $d$ pairs of
vertices, i.e. the graph with $V=\{1,2, \dots, 2d\}$ and where $E$ contains all
possible edges except those of the form $\{i,i+1\}$ for $i$ odd (see Figure~\ref{fig:K}). 
The graph $K_d(2)$ is $(2,3)$-agreeable, and by \cite{roberts}*{Theorem~7}, we have $\bx(K_d(2))=d=\#V/2$.
\end{example}

\begin{rem}
Graphs with $\bx(G)\leq 1$ are \emph{interval graphs}, which can
be easily identified in linear time~\cites{interval,interval2}. Algorithms exist to test if $\bx(G) \leq 2$~\cite{dleq2}, or to compute boxicity in general~\cite{compute}, but they are a lot more cumbersome. The task of testing if
$\bx(G)\leq d$ is known to be NP-complete for all $d\geq 2$~\cite{compute}.
\end{rem}

The definition of $(2,3)$-agreeability as it appears in~\cite{agree} can be reformulated in terms of intersection graphs.
\begin{definition}\label{df:23}
An arrangement $\B=\{B_1, \dots, B_n\}$ of $d$-boxes
is $(2,3)$-agreeable
if and only if any one of the three equivalent properties holds: 
\begin{enumerate}%
\item For any $1 \leq i<j<k \leq n$, one at least of the 
intersections $B_i\cap B_j$, $B_i\cap B_k$ or $B_j\cap B_k$ 
is non-empty.
\item For any three vertices in the intersection graph $G_\B$, the graph induced by these vertices contains at least one edge.
\item The graph complement $\overline{G_\B}$ of the intersection 
graph verifies $\omega(\overline{G_\B}) \leq 2$.
\end{enumerate}
\end{definition}

\subsection{Agreement Number and Agreement Proportion}
Since any simple, undirected graph can be realized as 
the intersection graph of an arrangement of boxes, 
it will be convenient to blur the distinction between 
the two notions. In particular, we can use 
properties~(2) and~(3) in Definition~\ref{df:23}
to define $(2,3)$-agreeability for graphs rather than arrangement.

Another good reason to identify arrangements and their graphs
is that the intersection graph encodes all the information 
about arrangements of boxes (this fails for arrangements of more general convex sets). Indeed, in such an arrangement, having nonempty pairwise intersection and having a point common to all the boxes are equivalent. In particular, the maximal number of overlapping boxes (or \emph{agreement number} of the society) is simply the clique number 
$\omega(G_\B)$ of the intersection graph. 
%% not used!!
%To improve readability, 
%we use the notation $\omega(\B)=\omega(G_\B)$.

\begin{notation}
We denote by $\G$ the set of all $(2,3)$-agreeable graphs, and,  
for any $d\geq 0$, denote by $\G_d$ the subset of those graphs with boxicity at most $d$. Given $r\geq 1$, we let $\G(r)$ and $\G_d(r)$ 
respectively be the subsets of $\G$ and $\G_d$ 
formed by graphs $G$ with 
$\omega(G)\leq r$.
Note that for any $G\in\G_d(r)$ and any subset of vertices $W\sub V(G)$, the subgraph $G[W]$
induced by $W$ is also in $\G_d(r)$: $(2,3)$-agreeability is preserved by taking induced graphs, and both clique size and boxicty
can only decrease (Remark~\ref{rem:subgraph}).
\end{notation}

We define the associated vertex sizes for all $r \geq 1$ and all $d \geq 0$,
\begin{align*}
\eta(r,d) &= \max\{\#V(G) \mid G \in \G_d(r)\}; \\
\eta(r) &= \max\{\#V(G) \mid G \in \G(r)\}.
\end{align*}
These quantities are related by the inequalities
$$ 2r=\eta(r,1)\leq \eta(r,2) \leq \dots \leq \eta(r).$$
We will show in Proposition~\ref{prop:growth_eta} that $\eta(r)$
is finite for all $r\geq 1$, and thus that all sets $\G_d(r)$ are finite too. This is not a surprising result, since it is the 
expected behavior brought on by $(k,m)$-agreeability; but note that, in our case of interest, the very existence of a positive agreement proportion was left open in~\cite{agree}.% (see \S~1).
\bigskip

For any graph, the \emph{agreement proportion} is defined as 
$\omega(G)/\#V(G).$ Once we prove that the set $\G_d(r)$ 
is finite for all $r\geq 1$ and $d\geq 1$, we can define
\begin{equation}\label{eq:rho}
\rho(r,d)=\min\left\{\omega(G)/\#V(G) \mid G\in \G_d(r)\right\},
\end{equation}
i.e. the \emph{minimal agreement proportion} that can be obtained 
from a $(2,3)$-agreeable graph
with boxicity at most $d$ and clique number at most $r$.

%%%%%%%%%%%%%%%%%%%%%%%%%%%%%%%%%%%%%%%%%%%%%%%%%%%%%%%%%%%%%%%%%%

\section{Upper and Lower Bounds on Degrees}

Throughout this section,  $G=(V,E)$ denotes a $(2,3)$-agreeable graph
on $n$ vertices. We show that a $(2,3)$-agreeable graph with low clique number must have many edges.
The results obtained here are purely combinatorial: in this section, we ignore the geometry of the problem  
and the boxicity of $G$.

\subsection{Lower Bound on the Degree}

The following trivial observation is the key to establishing lower bounds on the degrees of vertices.

\begin{lemma}\label{lem:deg}
If $G$ is a $(2,3)$-agreeable graph, then for 
any vertex $v \in V$, we have $\deg(v)
\geq n-\omega(G)-1.$
\end{lemma}
Note that the inequality in this lemma may be strict, even if $v$ is of minimal degree. We can see this by considering $G=W_4$, the wheel with four spokes, which is a $(2,3)$-agreeable graph with $n=5$ and 
$\omega(G)=3$.

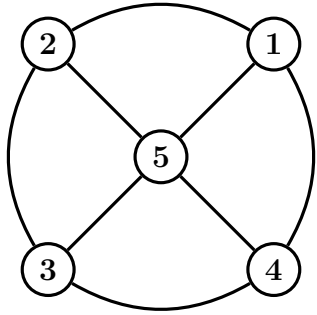
\begin{figure}
\centering
\begin{tikzpicture}[scale=1.5, style=very thick]
\tikzstyle{every node}=[draw,shape=circle];
\path (1,1) node  (a) {\lsty{1}};
\path (-1,1) node  (b) {\lsty{2}};
\path (-1,-1) node  (c) {\lsty{3}};
\path (1,-1) node   (d) {\lsty{4}};
\path (0,0) node   (t) {\lsty{5}};
\draw (a) to[bend right] (b);
\draw (b) to[bend right] (c);
\draw (c) to[bend right] (d);
\draw (d) to[bend right] (a);
\draw (a) -- (t) -- (c);
\draw (b) -- (t) -- (d);
\end{tikzpicture}
\caption{The wheel with four spokes $W_4$ is an example of a graph for 
which the inequality in Lemma~\ref{lem:deg} is strict.}
\end{figure}

\begin{proof}
The vertex $v \in V$ is connected to $\deg(v)$ vertices.  The
remaining $n-\deg(v)-1$ vertices must form a clique $W$. 
Indeed, if $W$ was not a clique, it would contain two non-adjacent vertices, $u$ and $w$. The subgraph induced by the three
vertices $\{u,v,w\}$ would be empty, which would contradict the fact that $G$ is  $(2,3)$-agreeable.
Thus, $\omega(G)\geq |W|=n-\deg(v)-1$, and the result follows.
\end{proof}

Using the formula 
$$ |E|=\frac{1}{2} \sum_{v\in V} \deg(v),$$
Lemma~\ref{lem:deg} 
yields the following lower-bound on $|E|$.
\begin{corollary}\label{cor:lower}
For any $(2,3)$-agreeable graph $G$, we have
$$|E| \geq \frac{n}{2}(n-\omega(G)-1).$$
\end{corollary}
The 5-cycle from Remark~\ref{rem:Z5} verifies
$n=5$, $\omega(G)=2$ and $|E|=5$.
This shows that the bound in Corollary~\ref{cor:lower} is sharp. 

\subsection{Examples with Low Agreement Proportion} 
The conclusion of~\cite{agree} mentioned the existence of 
$(2,3)$-agreeable families of $2$-boxes with agreement $3/8$.
(The example, credited to Rajneesh Hegde, was not given in the paper.)
We give a few examples.

\begin{figure}
\centering
\begin{tikzpicture}[rotate=90, scale=1.5, style=very thick]
\tikzstyle{every node}=[draw,shape=circle];
\path (0,0) node [fill=violet, fill opacity=1] (a) {\lsty{3}};
\path (0,1) node [fill=blue, fill opacity=1] (b) {\lsty{2}};
\path (1,1) node [fill=orange, fill opacity=1] (c) {\lsty{4}};
\path (1,0) node  [fill=purple, fill opacity=1] (d) {\lsty{1}};
\path (-1,-1) node [fill=green, fill opacity=1]  (t) {\lsty{7}};
\path (-1,2) node [fill=brown, fill opacity=1] (u) {\lsty{8}};
\path (2,2) node [fill=gray, fill opacity=1] (v) {\lsty{6}};
\path (2,-1) node [fill=yellow, fill opacity=1] (w) {\lsty{5}};
\draw (a) -- (b) -- (c) -- (d) -- (a);
\draw (t) -- (u) -- (v) -- (w) -- (t);
\draw (v) -- (c) -- (u) -- (b) -- (t) -- (a) -- (w) -- (d) -- (v);
\end{tikzpicture}
\hfil
\begin{tikzpicture}[scale=.25, style=very thick]
\draw[fill=purple, fill opacity=0.8] (5,0) rectangle (6,20); %1
\draw[fill=blue, fill opacity=0.8] (15,0) rectangle (16,20); %2
\draw[fill=violet, fill opacity=0.8] (0,5) rectangle (20,6); %3 
\draw[fill=orange, fill opacity=0.8] (0,15) rectangle (20,16); %4
\draw[fill=yellow, fill opacity=0.3] (2,1) rectangle (14,10); %5
\draw[fill=gray, fill opacity=0.3] (1,7) rectangle (8,18); %%% 6
\draw[fill=green, fill opacity=0.2] (12,4) rectangle (17,14); %7
\draw[fill=brown, fill opacity=0.2] (7,13) rectangle (19,19); %8
\end{tikzpicture}
\caption{A $(2,3)$-agreeable graph with $|V(G)|=8$, $\omega(G)=3$ and $\bx(G)=2$, together with a family of $2$-boxes whose intersection graph is $G$. This graph is $4$-regular, $|E|=16$.}\label{fig:38a}
\end{figure}
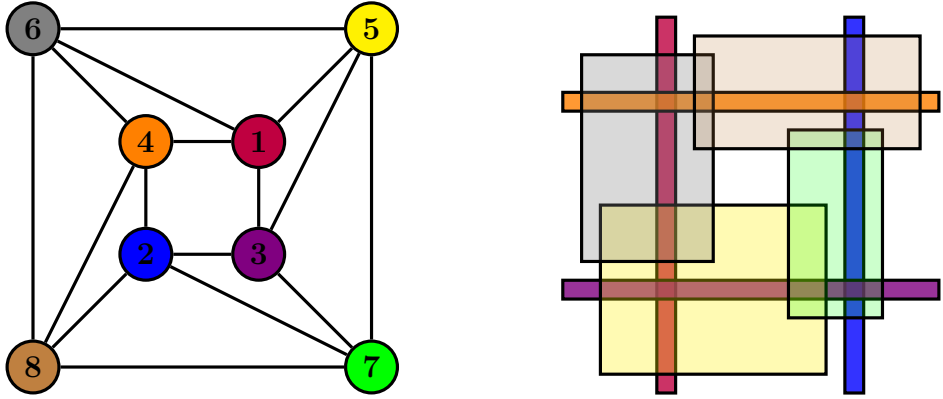

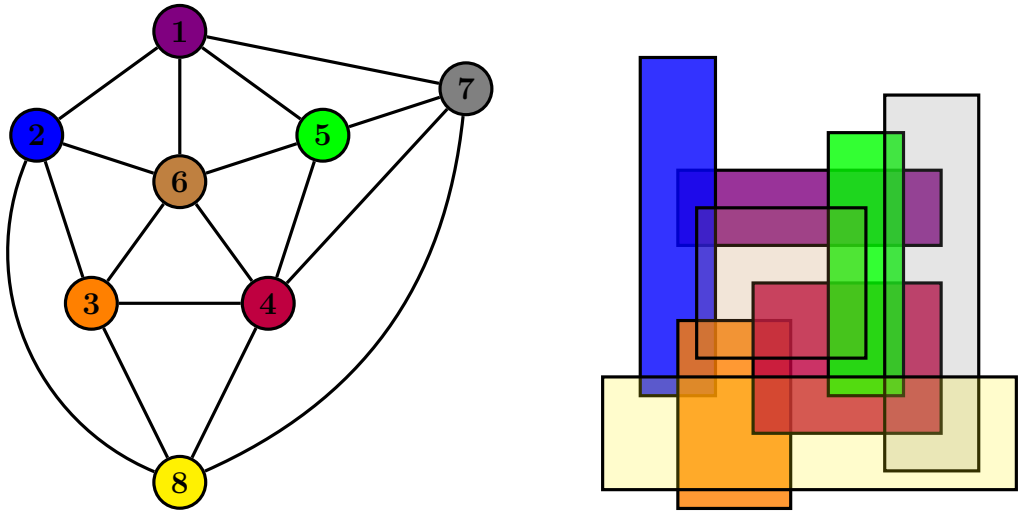
\begin{figure}
\centering
\begin{tikzpicture}[rotate=90, scale=2, style=very thick]
\tikzstyle{every node}=[draw,shape=circle];
\path (0:0cm) node [fill=brown, fill opacity=1] (u) {\lsty{6}};
\path (0:1cm) node [fill=violet, fill opacity=1] (a) {\lsty{1}};
\path (72:1cm) node [fill=blue, fill opacity=1] (b) {\lsty{2}};
\path (2*72:1cm) node [fill=orange, fill opacity=1] (c) {\lsty{3}};
\path (3*72:1cm) node  [fill=purple, fill opacity=1] (d) {\lsty{4}};
\path (4*72:1cm) node [fill=green, fill opacity=1]  (e) {\lsty{5}};
\path (288:2cm) node [fill=gray, fill opacity=1] (v) {\lsty{7}};
\path (180:2cm) node [fill=yellow, fill opacity=1] (w) {\lsty{8}};
\draw (a) -- (b) -- (c) -- (d) -- (e) -- (a);
\draw (a) -- (u); \draw (b) -- (u);
\draw (c) -- (u); \draw (d) -- (u); \draw (e) -- (u);
\draw (a) -- (v); \draw (e) -- (v); \draw (d) -- (v);
\draw (w) -- (c); \draw (w) -- (d);
\draw (w) to [bend left=45] (b);
\draw (w) to [bend right] (v);
\end{tikzpicture}
\hfil
\begin{tikzpicture}[scale=.5, style=very thick]
\draw[fill=violet, fill opacity=0.8] (1,7) rectangle (8,9); %1
\draw[fill=blue, fill opacity=0.8] (0,3) rectangle (2,12); %2
\draw[fill=orange, fill opacity=0.8] (1,0) rectangle (4,5); %3 
\draw[fill=purple, fill opacity=0.8] (3,2) rectangle (8,6); %4
\draw[fill=green, fill opacity=0.8] (5,3) rectangle (7,10); %5
\draw[fill=brown, fill opacity=0.2] (1.5,4) rectangle (6,8); %%% 6
\draw[fill=gray, fill opacity=0.2] (6.5,1) rectangle (9,11); %7
\draw[fill=yellow, fill opacity=0.2] (-1,0.5) rectangle (10,3.5); %8
\end{tikzpicture}
\caption{Another $(2,3)$-agreeable graph with $|V(G)|=8$, 
$\omega(G)=3$ and $\bx(G)=2$. Like 
the example from Figure~\ref{fig:38a}, t
his graph also has boxicity~2, but $|E|=17$.
}
\label{fig:38b}
\end{figure}

\begin{figure} \centering
\begin{tikzpicture}[scale=1.5, style=very thick]
\tikzstyle{every node}=[draw,shape=circle, minimum size = 25pt];
\path (135:4cm) node (1) {\lsty{1}};
\path (45:4cm) node (2) {\lsty{2}};
\path (225:4cm) node (3) {\lsty{3}};
\path (315:4cm) node (4) {\lsty{4}};
\path (180:2cm) node (5) {\lsty{5}};
\path (90:2cm) node (6) {\lsty{6}};
\path (0:2cm) node (7) {\lsty{7}};
\path (270:2cm) node (8) {\lsty{8}};
\path (135:1.414cm) node (9) {\lsty{9}};
\path (45:1.414cm) node (10) {\lsty{10}};
\path (315:1.414cm) node (11) {\lsty{11}};
\path (225:1.414cm) node (12) {\lsty{12}};
\path (0:0cm) node (13) {\lsty{13}};
\draw[color=green](1) -- (2); 
\draw[color=green](1) -- (3); 
\draw[dashed, color=green](1) -- (7); 
\draw[color=green](1) -- (9); 
\draw[color=green](2) -- (4); 
\draw[dashed, color=green](2) -- (8); 
\draw[color=green](2) -- (10); 
\draw[color=green](3) -- (4); 
\draw[dashed, color=green](3) -- (6); 
\draw[color=green](3) -- (12); 
\draw[dashed, color=green](4) -- (5); 
\draw[color=green](4) -- (11); 
\draw[dashed, color=green](5) to[bend right=25] (7); 
\draw[color=green](5) -- (9); 
\draw[color=green](5) -- (12); 
\draw[dashed, color=green](6) to[bend left=25] (8); 
\draw[color=green](6) -- (9); 
\draw[color=green](6) -- (10); 
\draw[color=green](7) -- (10); 
\draw[color=green](7) -- (11); 
\draw[color=green](8) -- (11); 
\draw[color=green](8) -- (12); 
\draw[color=green](9) -- (13); 
\draw[color=green](10) -- (13); 
\draw[color=green](11) -- (13); 
\draw[color=green](12) -- (13);
\end{tikzpicture}

\caption{The \emph{complement} of this graph is a $(2,3)$-agreeable $8$-regular graph with $\omega=4$.
}
\label{fig:134}
\end{figure}
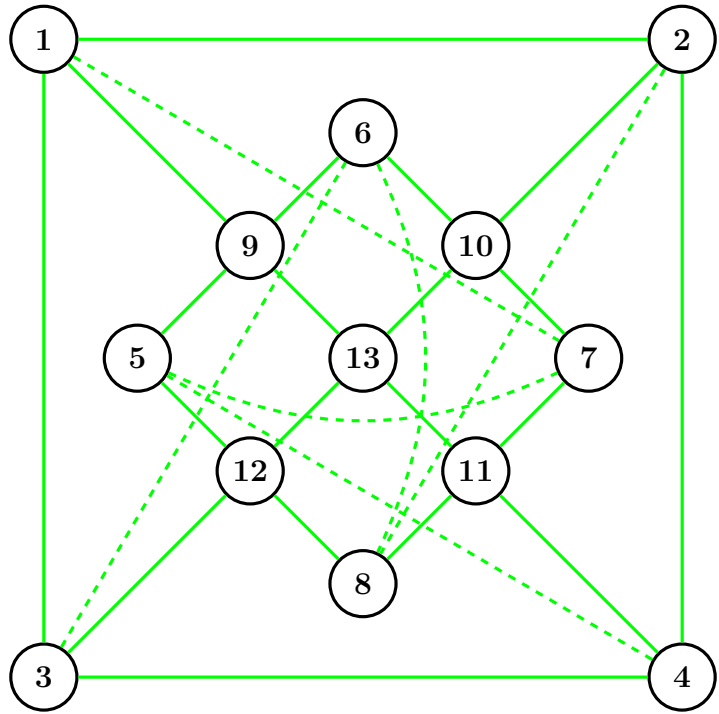

\begin{description}
\item[Case $n=8, \bx(G)=2$.] 
Figures~\ref{fig:38a} and~\ref{fig:38b}
give two non-isomorphic examples of families of eight $2$-boxes with 
no more than triple intersections. The corresponding intersection graphs have respectively $8$ and $10$ triangles.
\item[Case $n=8, \bx(G)=?$] Figure~\ref{fig:38c} presents a third example of a $(2,3)$-agreeable graph with agreement proportion $3/8$, obtained from Figure~\ref{fig:38a} by adding two edges. 
This graph has 12 triangles;
its boxicity may be more than 2 (we conjecture that it is).
\item[Case $n=13, \omega(G)=4$.]
Figure~\ref{fig:134} presents (the complement of) a $(2,3)$-agreeable graph on $13$ vertices with unknown boxicity and agreement proportion $4/13\approx 0.31$. There are $39$ distinct cliques of size~4 in that example. We will prove in Proposition~\ref{prop:table} that no $(2,3)$-agreeable graph on 14 or more vertices has such low clique number.
\end{description}

\begin{figure}
\centering
\begin{tikzpicture}[rotate=90, scale=1.5, style=very thick]
\tikzstyle{every node}=[draw,shape=circle];
\path (0,0) node [fill=violet, fill opacity=1] (a) {\lsty{3}};
\path (0,1) node [fill=blue, fill opacity=1] (b) {\lsty{2}};
\path (1,1) node [fill=orange, fill opacity=1] (c) {\lsty{4}};
\path (1,0) node  [fill=purple, fill opacity=1] (d) {\lsty{1}};
\path (-1,-1) node [fill=green, fill opacity=1]  (t) {\lsty{7}};
\path (-1,2) node [fill=brown, fill opacity=1] (u) {\lsty{8}};
\path (2,2) node [fill=gray, fill opacity=1] (v) {\lsty{6}};
\path (2,-1) node [fill=yellow, fill opacity=1] (w) {\lsty{5}};
\draw (a) -- (b) -- (c) -- (d) -- (a);
\draw (t) -- (u) -- (v) -- (w) -- (t);
\draw (v) -- (c) -- (u) -- (b) -- (t) -- (a) -- (w) -- (d) -- (v);
\draw (b) -- (d);
%\draw[dashed] (v) to[bend left=135] (t);
\draw (v) to[bend left=55] (3,-2) to[bend left=55] (t);
\end{tikzpicture}
\caption{A $(2,3)$-agreeable graph with $|V(G)|=8$, $\omega(G)=3$,
but $\bx(G)$ unknown.
Modifying the arrangement of Figure~\ref{fig:38a} to have 
$B_1 \cap B_2\neq \varnothing$ and $B_6 \cap B_7\neq \varnothing$
creates more intersections, so it is not immediately obvious whether 
a $2$-box arrangement can realize this graph.
 }\label{fig:38c}
\end{figure}
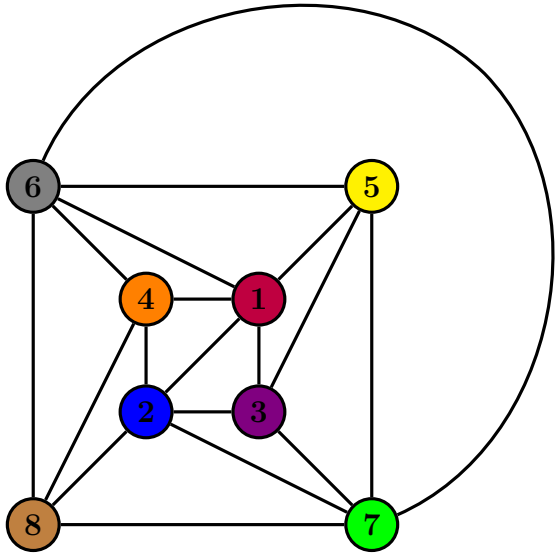

\subsection{Upper Bounds on Degree and Graph Size}
We now give upper-bounds on the degrees of vertices in 
$(2,3)$-agreeable graphs, and deduce an upper bound on 
the number of vertices of such a graph with given clique number.

\begin{lemma}\label{lem:Delta}
Let $G\in\G_d(r)$, where $r\geq 2$ and $d\geq 1$. Then 
for any $v\in V$, we have 
$$\deg(v) \leq \eta(r-1,d) \leq \eta(r-1).$$
\end{lemma}

\begin{proof}
The neighbors of $v$ induce a $(2,3)$-agreeable graph $H$. If there are more than $\eta(r-1,d)$ vertices in the graph $H$, it must contain an $r$-clique, which together with $v$ forms an $(r+1)$-clique in $G$, contradicting the hypothesis $\omega(G)\leq r$.
\end{proof}
The inequality in the lemma can be sharp, but it is not always so, 
even if $G$ has the maximum $\eta(r)$ vertices: taking $r=3$, 
we'll see in Proposition~\ref{prop:table} that $\eta(r-1)=5$ 
and $\eta(r)=8$. The graphs in
Figure~\ref{fig:38a} and in Figure~\ref{fig:38b} 
both have the maximum 8~vertices for their clique number of~3,
but the maximum degree is~4 in the first graph, but a sharp~5 in the 
second example.
%has 8~vertices but is only 4-regular, whereas the graph in Figure~\ref{fig:38b} has also 8~vertices, including one of degree~5.

\bigskip

With Lemma~\ref{lem:deg} giving a lower bound on the number of edges which increases with the number of vertices, and Lemma~\ref{lem:Delta} giving an upper-bound which depends only on the clique number, 
this suggests that the graphs which maximize $\eta(r)$ must 
be regular or almost-regular.
We can use this idea to, step by step, establish the first few values of $\eta(r)$. 

\begin{proposition}\label{prop:table}
We have the following table for the maximal size of $(2,3)$-agreeable graphs with $\omega(G)=r$.
\begin{center}
\begin{tabular}{c||*{5}{c|}}
$r$ & $1$ & $2$ & $3$ & $4$ & $5$ \\
\hline
$\eta(r)$ & $2$ & $5$ & $8$ & $13$ & $\leq 18$ \\
\end{tabular}
\end{center}
\end{proposition}

\begin{proof}
By definition of $(2,3)$-agreeability, any graph with at least three vertices must have an edge, and thus $\eta(1)=2$. The examples we've seen so far give the following lower bounds:
$$ \eta(2) \geq 5, \quad
\eta(3) \geq 8, \quad
\eta(4) \geq 13.
$$ 
Suppose that these lower bounds are not sharp: i.e. there exists $(2,3)$-agreeable graphs with the following.
\begin{center}
\begin{tabular}{c||*{4}{c|}}
$\#V$ &  $2$ & $3$ & $4$ & $5$ \\
\hline
$\omega(G)$ & $6$ & $9$ & $14$ & $19$ \\
\end{tabular}
\end{center}
Let $\delta(G)$  and $\Delta(G)$ denote respectively the minimum and the maximum degree for vertices in $G$. The case $|V(G)|=6$ and  $\omega(G)=2$ is clearly impossible, since  Lemma~\ref{lem:deg} implies 
$\delta(G)\geq 3$, and Lemma~\ref{lem:Delta} implies 
$\Delta(G) \leq \eta(1)=2$, giving the contradiction 
$\delta(G)>\Delta(G)$.

Thus, we have proved that $\eta(2)=5$, which, combined with 
Lemma~\ref{lem:Delta} implies that for any $G\in \G$ with 
$\omega(G)=3$, we must have $\Delta(G)\leq 5$. In the case $|V(G)|=9$
with $\omega(G)=3$, 
Lemma~\ref{lem:deg} yields $\delta(G) \geq 5$. Since the graph $G$ cannot be $5$-regular (the sum of all degrees must be even), this yields in turn $\Delta(G)\geq 6$ which is again a contradiction.

This proves $\eta(3)=8$, which implies that $\Delta(G)\leq 8$ for any $G\in \G$ with $\omega(G)=4$. The other cases are similar.
\end{proof}

The method used in the proof of the above proposition could be
extended indefinitely, provided one can construct examples that 
provide lower bounds on $\eta$.
Even without a battery of examples, we can  prove that  
the function $\eta(r)$ has at most quadratic growth.
Thus, the sets $\G_d(r)$ are finite for any $d\geq 1$ 
and $r \geq 1$.

\begin{proposition}\label{prop:growth_eta}
For all $r \geq 1$, the maximal number of vertices $\eta(r)$ 
for a $(2,3)$-agreeable graph $G$ with $\omega(G)\leq r$ 
verifies $\eta(r) \leq r(r+3)/2$.
\end{proposition}

\begin{proof}
Let $G$ be a $(2,3)$-agreeable graph such that $\omega(G)=r$ 
and $|V(G)|=\eta(r)$. If $v$ is a vertex of $G$, Lemma~\ref{lem:deg} and Lemma~\ref{lem:Delta} imply the inequalities:
$$ \eta(r)-r-1 \leq \deg(v) \leq \eta(r-1).$$
Solving the recurrence 
$ \eta(r)-r-1 - \eta(r-1)\leq 0$
with the initial condition $\eta(1)=2$ gives the result.
\end{proof}

\section{Lower bound on Boxicity and the Main Result}

Given a simple graph $G$ on $n$ vertices, call a vertex $v\in V(G)$ 
\emph{universal} if $\deg(v)=n-1$. 
The preprint~\cite{lower} presents several lower bounds on 
the boxicity of a graph; we will need the following one.

\begin{theorem}[\cite{lower}]\label{thm:lower}
Let $G$ be a graph with no universal vertices and minimum degree 
$\delta$. Then the boxicity of $G$ has the lower bound:
$$ \bx(G) \geq \frac{n}{2(n-\delta-1)}.$$
\end{theorem}
The theorem above only applies to graphs with no universal vertices.
Fortunately, the lemma below shows we only need 
to consider such graphs.
Recall that for all $r\geq 1$ and $d \geq 1$, the quantity $\rho(r,d)$ denotes the minimum agreement proportion that can be achieved by a graph $G\in\G_d(r)$.

\begin{lemma}\label{lem:rho}
Given $r\geq 1$ and $d \geq 1$, consider a graph $G\in \G_d(r)$ 
such that the agreement proportion of $G$ is equal to $\rho(r,d)$. Then, $G$ has no universal vertices.
\end{lemma}

\begin{proof}
Suppose $G\in \G_d(r)$ is a graph with universal vertices, 
$G\neq K_n$.
We construct from $G$ a graph $\GG\in \G_d(r)$ without universal vertices and with a lower agreement proportion.
Let $\Omega$ be the set of universal vertices,
$$\Omega=\{v\in V(G)\mid \deg(v)=n-1\};$$ 
define $W=V(G)\bs \Omega$ and let $\GG=G[W]$ be the graph induced by $W$.
Since we assumed $G\neq K_n$, the graph $\GG$ is non-empty.
Note that $\bx(\GG)\leq \bx(G)\leq d$, since boxicity can only decrease when considering induced graphs (Remark~\ref{rem:subgraph}).
Letting $k=|\Omega|$, we have for any vertex in $w\in W$, 
$$ \deg_{\GG}(w)=\deg_G(w)-k<n-1-k=|W|-1,$$
so that no vertex in $\GG$ is universal. Moreover, we have 
$$\omega(\GG)=\omega(G)-k,$$
since any maximal clique in $G$ must contain all the vertices in 
$\Omega$.
Thus, the agreement proportion for $\GG$ is
$$ \frac{\omega(\GG)}{\#V(\GG)}=\frac{\omega(G)-k}{n-k}< \frac{\omega(G)}{n}; $$
thus, any graph which minimizes agreement proportion does
not have any universal vertices. 
\end{proof}

Our main result, Theorem~\ref{thm:main}, is now easily derived.

\begin{theorem}
For any $r\geq 1$ and $d\geq 1$, we have
$$ \rho(r,d)\geq \frac{1}{2d}.$$
In particular, any $(2,3)$-agreeable arrangement of $d$-boxes 
must have an agreement proportion of at least~$(2d)^{-1}$.
\end{theorem}

\begin{proof}
Consider a graph $G\in \G_d(r)$ on $n$ vertices such that the agreement proportion of $G$ is equal to the minimum $\rho(r,d)$.
By Lemma~\ref{lem:rho}, $G$ does not contain any universal vertex.
Theorem~\ref{thm:lower} applies, so that 
$$ d \geq \bx(G) \geq \frac{n}{2(n-\delta-1)};$$
where $\delta$ denotes the minimum degree in~$G$.
Since $G$ is $(2,3)$-agreeable, Lemma~\ref{lem:deg} yields 
$$\omega(G) \geq n-\delta-1.$$
Combining the two inequalities, we get
$$ \rho(r,d)=\frac{\omega(G)}{n}
\geq \frac{n-\delta-1}{n}\geq \frac{1}{2d}.$$
This completes the proof of the main theorem.
\end{proof}

\section{Some Questions}
Our research did not yield any general method to construct $(2,3)$-agreeable graphs with low agreement numbers. 
Achieving this while keeping a handle on boxicity is even more of a challenge, especially given the hardness of computations.
As noted in~\cite{lower}, upper bounds on $\bx(G)$ have been extensively studied, but results about lower bounds are scarcer, 
and any new development in this direction could conceivably impact this work.

\medskip

The sharpness of our bounds for $d\geq 2$ remains unknown. 
Answers to the following questions would have a great impact on 
the total understanding of the $(2,3)$-agreeable case.

\begin{description}
\item[Are there examples with high boxicity?]
Example~\ref{ex:K} shows that it is easy to find $(2,3)$-agreeable graphs with arbitrarily high boxicity. But the $d$-partite graph $K_d(2)$ has a high agreement proportion, $1/2$. It might be that the twin constraints of $(2,3)$-agreeability and low clique number are somehow at odds with having high boxicity. It might be that the sequence of maximal sizes of  $(2,3)$-agreeable graphs 
with $\bx(G)\leq d$ and $\omega(G)\leq r$,
$$ 2r=\eta(r,1) \leq \eta(r,2) \leq \eta(r,3) \leq \dots $$
becomes eventually constant for $d \geq d_0$, for a value $d_0$ \emph{independent of $r$}. (When $r$ is fixed, we are dealing with a finite number of graphs, so the  boxicity is trivially bounded.)
If this sequence indeed stabilizes, it implies that $(2,3)$-agreeable arrangement of boxes have a positive agreement proportion which is independent of $d$.

\item[Does the agreement proportion go to zero with $d$?]
The lower bound on the agreement proportion of $(2d)^{-1}$, obtained in Theorem~\ref{thm:main}, goes to zero when $d$ goes to infinity. 
But even if the true agreement proportion is a strictly decreasing function of $d$, it might still have a positive limit. In that case, we would again have a positive lower bound valid for all $d$.

If it were the case, it would imply that
$\eta(r)$ does \emph{not} grow quadratically, 
as Proposition~\ref{prop:growth_eta} suggests it could.

\item[Graph agreement proportion.] Another way to ask the same questions is the following: let 
$$ \rho(r) = \rho\left(r, \lfloor \eta(r)/2 \rfloor \right). $$
%\min_{1 \leq d \leq \lfloor \eta(r)/2 \rfloor}\rho(r,d).$$
Since any $(2,3)$-agreeable graph has at most $\eta(r)$ vertices 
and boxicity as most half its number of vertices, 
$\rho(r)$ is the well-defined minimum agreement proportion 
obtainable from a $(2,3)$-agreeable graph $G$ with 
$\omega(G) \leq r$. The problem becomes to understand how $\rho(r)$ 
varies as a function of $r$.
\end{description}

\appendix 

\section{The Exposed Box Method}
During the undergraduate research project when this research was started, in summer 2008, we had originally obtained a different  bound, using another method. The results are comparable for low values of the boxicity, but Theorem~\ref{thm:main} always gives a stronger result, and the difference quickly becomes large (see Figure~\ref{fig:table}).  
The method was adapted 
from the paper~\cite{boxes1} by Eckhoff, which deals with maximizing the entries in the face vector of an arrangement of $d$-boxes.
Since that method is of independent interest, and may be extended in other contexts, we outline the result and the main steps of
 the proof in this appendix.
(The reader will notice that Corollary~\ref{cor:lower} is the only required result which is specific to $(2,3)$-agreeability. Every tool presented in the appendix can be extended to the $(k,m)$-agreeable case.)

\begin{theorem}\label{thm:oldmain}
For any $d \geq 1$, there exists  $\betabox(d)>0$ such that 
any $(2,3)$-agreeable $d$-box society of size $n$ has an overlap of size at least $\betabox(d)n$. Moreover, 
$\betabox(d) \geq F^{[d-1]}(1/2),$ 
where $F^{[d-1]}$ denotes the function
\begin{equation*}
F(x)=\frac{-2 - x + \sqrt{4 - 4 x + 5 x^2}}{2 (x-2)},
\end{equation*}
iterated $d-1$ times. In particular, for $(2,3)$-agreeable boxes in 
$\R^2$, the agreement proportion verifies 
$$\betabox(2) \geq 
\frac{5-\sqrt{13}}{6} \approx 0.2324.$$
\end{theorem}

\begin{figure}[ht]
\renewcommand{\arraystretch}{1.2}
\centering
\begin{tabular}{c||*{5}{c|}}
$d$ & 1 & 2 & 3 & 4 & 5 \\
\hline
$(2d)^{-1}$ & 0.5 & 0.25 & 0.167 & 0.125 & 0.1 \\
\hline
$F^{[d-1]}(1/2)$ & 0.5 &0.23 & 0.11 & 0.05 & 0.02\\
\hline
\end{tabular}
\caption{A comparison of the results given by Theorem~\ref{thm:main} and Theorem~\ref{thm:oldmain}, for boxicity up to~5.}
\label{fig:table}
\end{figure}

\subsection{The Eckhoff Induction}
Let $\B=\{B_1, \dots, B_n\}$ denote an
arrangement of $n$ boxes in $\R^d$.  
Recall that if $Q \subseteq \R^d$ is a convex subset, a hyperplane $H$ is said to be a \emph{supporting hyperplane} for $Q$ if $Q\cap H \neq \varnothing$ and $Q$ is entirely contained in one of the two half-spaces delimited by $H$.
We borrow the
following key notion from~\cite{boxes1}. 

\begin{definition}
We say that the box $Q\in \B$ is \emph{exposed} by the hyperplane $H$
for the arrangement $\B$ if 
\begin{enumerate}\itemsep=0pt
\item $H$ is a \emph{supporting hyperplane} of $Q$ which is parallel to some
coordinate hyperplane. 
\item For any $P\in \B$ such that $P\cap H=\varnothing$, $P$ and $Q$
do not lie in the same half-space defined by $H$. 
\end{enumerate}
We call $H$ an \emph{exposing hyperplane} for $Q$.
\end{definition}

\begin{figure}[ht]
\centering

\begin{tikzpicture}[scale=1.25]
\draw[fill=orange, fill opacity=0.3] (0.5,0.5) rectangle (1.25,1.5); %E
\draw[fill=blue, fill opacity=0.3] (1,0) rectangle (2,3); %F
\draw[fill=purple, fill opacity=0.3] (0,1.25) rectangle (4,2); %D
\draw[fill=green, fill opacity=0.3] (.5,2.5) rectangle (1.75,3.5);%A
\draw[fill=gray, fill opacity=0.3] (2.25,1.75) rectangle (4.25,4);%B
\draw[fill=yellow, fill opacity=0.3] (2.5,.5) rectangle (3.5,2.2);%C
\draw[style=very thick] (-1,2.5) -- (5,2.5);
\path (5.25,2.5) node (h) {$H$};
\path (.75,3.25) node (a)  {$A$};
\path (2.5,3.75) node (b)  {$B$};
\path (3.25,.75) node (c)  {$C$};
\path (.25,1.75) node (d)  {$D$};
\path (.75,.75) node (e)  {$E$};
\path (1.75,.25) node (f)  {$F$};
\end{tikzpicture}

\caption{The Box $A$ is exposed by $H$. The boxes $B$ and $F$ 
meet $H$, and the other three are in the lower half-plane.}
\end{figure}
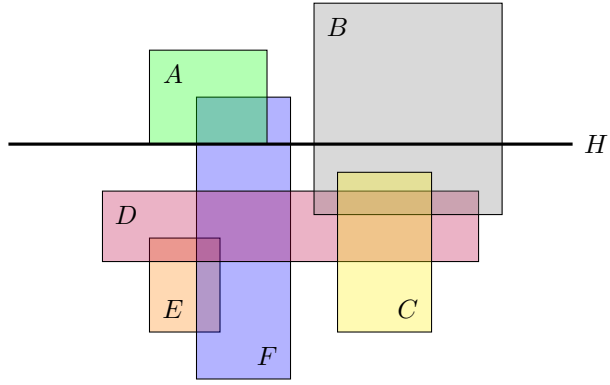

Note that any arrangement of boxes $\B$ contains exposed boxes: 
take any coordinate hyperplane $H$ infinitely far from $\B$, then bring it back towards the 
arrangement. The first box which $H$ supports is exposed by that hyperplane. If $B_1$ is exposed in $\B$ by some hyperplane $H$, 
we let 
\begin{equation}\label{eq:BB}
\B'=\{B_2, \dots,B_n\} \quad \text{and}\quad 
\B''=\{B_1 \cap B_i \mid 2 \leq i \leq n\}.
\end{equation}
For any arrangement $\B$, we denote by $f_k(\B)$ 
the number of non-empty 
$(k+1)$-fold intersections in $\B$.
In particular, $f_0(\B)$ is the number of boxes in the arrangement
and $f_1(\B)$ is the number of pairwise intersection, i.e. the 
number of edges in the intersection graph $G_\B$.

\begin{lemma}\label{lem:BBB}
If $B_1$ is exposed in $\B$, then for 
all $1 \leq k \leq n-1$, we have
\begin{equation*}
f_k(\B)=f_k(\B')+f_{k-1}(\B'').
\end{equation*}
\end{lemma}

The arrangements $\B'$ and $\B''$ being also $(2,3)$-agreeable,
the strategy is to apply Lemma~\ref{lem:BBB} in the case $k=1$, 
in order to relate number of vertices and number of edges. 
Denote by $e(n,r,d)$ the maximal number of edges for a graph $G\in\G_d(r)$ with $n$ vertices. We have the following recurrence relation:

\begin{lemma}\label{lem:e}
For all $n\geq r \geq 2$ and all  $d\geq1$,
we have
\begin{equation} e(n,r,d) \leq e(n-1,r,d)+\eta(r-1,d-1). \end{equation}
\end{lemma}

\subsection{An Induction on the Number of Boxes}

The recurrence on $n$ obtained in Lemma~\ref{lem:e} can be solved to 
relate the maximum number of edges $e(n,r,d)$ to the agreement proportion in the previous boxicity, $\betabox(d-1)$, and to obtain
the following.

\begin{theorem}\label{thm:e}
For all $n\geq r \geq 2$ and all  $d\geq1$, we have
\begin{equation}\label{eq:esolved}
e(n,r,d) \leq \binom{r}{2}+(n-r)\frac{r-1}{\betabox(d-1)}.
\end{equation}
\end{theorem}

\subsection{The Quadratic Formula and an Induction on the Boxicity}
Theorem~\ref{thm:oldmain} is proved by induction on $d$. The key is to compare the upper bound 
(Theorem~\ref{thm:e})
and the lower bound 
(Corollary~\ref{cor:lower}) on the number of edges to narrow down the possible values for $\omega(G)$.

\bigskip
\noindent
\textbf{A Sketch of the Proof of Theorem~\ref{thm:oldmain}. }
Let $G\in\G_d(r)$ and let $n$ be the number of vertices of $G$. The number of edges of $G$ verifies
\begin{equation}\label{eq:f1G}
\frac{n}{2}(n-r-1) \leq \#E(G) \leq \binom{r}{2}+(n-r)\frac{r-1}{\betabox(d-1)}.
\end{equation}
Let $\betabox= \betabox(d-1)$, and we can rewrite~\eqref{eq:f1G} as a quadratic inequality in $r$ 
with a negative coefficient for $r^2$,
\begin{equation}\label{eq:quadratic}
   (\betabox-2)r^2+(2n+\betabox n+2-\betabox)r 
  +(-\betabox n^2-2n+\betabox n)
   \geq 0.
\end{equation}
Thus $r$ is greater or equal to the smaller of the two roots of the quadratic expression. Dividing by $n$ and taking the limit as $n$ goes to infinity, most terms go to zero and we obtain the comparatively simpler 
$$ \lim_{n \to \infty} \frac{r}{n} \geq F(\betabox(d-1)); $$
with $F(x)$ defined as in Theorem~\ref{thm:oldmain}.
Thus an induction on $d$ gives an answer in terms of successive compositions of $F$.
\eop
\end{document}